\documentclass{amsproc}
\usepackage{amssymb,array,curves,calc}%
\newtheorem{Lemma}{Lemma}[section]\newcommand{\bel}{\begin{Lemma}}\newcommand{\eel}{\end{Lemma}}
\newtheorem{Proposition}[Lemma]{Proposition}\newcommand{\bprop}{\begin{Proposition}}\newcommand{\eprop}{\end{Proposition}}
\newtheorem{Theorem}[Lemma]{Theorem}\newcommand{\bthe}{\begin{Theorem}}\newcommand{\ethe}{\end{Theorem}}
\newcommand{\bpr}{{\bf Proof:~}}\newcommand{\epr}{$\blacksquare$\\}
\newtheorem{Remark}[Lemma]{Remark}\newcommand{\beR}{\begin{Remark}\rm}\newcommand{\eeR}{\end{Remark}}
\newtheorem{Definition}[Lemma]{Definition}\newcommand{\bed}{\begin{Definition}}\newcommand{\eed}{\end{Definition}}
\newtheorem{Example}[Lemma]{Example}\newcommand{\bex}{\begin{Example}\rm}\newcommand{\eex}{\end{Example}}
\newtheorem{Corollary}[Lemma]{Corollary}\newcommand{\bcor}{\begin{Corollary}\rm}\newcommand{\ecor}{\end{Corollary}}
\newtheorem{Fact}[Lemma]{Fact}\newcommand{\bfact}{\begin{Fact}\rm}\newcommand{\efact}{\end{Fact}}

\newcommand{\beq}{\begin{equation}}\newcommand{\eeq}{\end{equation}}
\newcommand{\bem}{\begin{displaymath}}\newcommand{\eem}{\end{displaymath}}
\newcommand{\beqa}{\begin{eqnarray}}\newcommand{\eeqa}{\end{eqnarray}}
\newcommand{\bee}{\begin{enumerate}}\newcommand{\eee}{\end{enumerate}}
\newcommand{\bei}{\begin{itemize}}\newcommand{\eei}{\end{itemize}}
\newcommand{\bet}{\begin{tabular}{cccccccc}}\newcommand{\eet}{\end{tabular}}
\newcommand{\bpm}{\begin{pmatrix}}\newcommand{\epm}{\end{pmatrix}}
\newcommand{\bM}{\begin{matrix}}\newcommand{\eM}{\end{matrix}}
\newcommand{\ber}{\begin{array}{l}}\newcommand{\eer}{\end{array}}

\newcommand{\tinyM}{\scriptstyle}
\newcommand{\tinyT}{\scriptsize}
\newcommand{\tinyA}{\tinyM\text{\tinyT}}

\newcommand{\cO}{{\mathcal{O}}}

\newcommand{\mC}{\mathbb{C}}

\newcommand{\mP}{\mathbb{P}}\newcommand{\mS}{\mathbb{S}}
\newcommand{\mR}{\mathbb{R}}

\newcommand{\mZ}{\mathbb{Z}}

\newcommand{\al}{\alpha}\newcommand{\be}{\beta}\newcommand{\de}{\delta}
\newcommand{\Ga}{\Gamma}
\newcommand{\ep}{\epsilon}\newcommand{\ka}{\kappa}\newcommand{\si}{\sigma}

\newcommand{\taues}{\tau^{es}}

\newcommand{\tC}{\tilde{C}}
\newcommand{\tS}{\tilde{S}}

\newcommand{\li}{~\\ $\bullet$ }\newcommand{\ls}{~\\ $\star$ }\newcommand{\into}{\stackrel{i}{\hookrightarrow}}
\newcommand{\ra}{\!\!\rightarrow\!\!}

\newcommand{\bin}[2]{{#1\choose{#2}}}  \newcommand{\ol}[1]{\overline{#1}}

  \newcommand{\omp}{ordinary multiple point}

\newcommand{\suml}{\sum\limits}\newcommand{\cupl}{\mathop\cup\limits}\newcommand{\oplusl}{\mathop\oplus\limits}

\newcommand{\mesh}[7]
{ \put(#1,#2){\vector(1,0){#6}}  \put(#1,#2){\vector(0,1){#7}}
\setcounter{tempx}{#3}  \addtocounter{tempx}{1}
\setcounter{tempy}{#4}  \addtocounter{tempy}{1}
\multiput(#1,#2)(#5,0){\value{tempx}}{\multiput(-1.5,-0.5)(0,#5){\value{tempy}}{.}}
}
\newcommand{\triang}[9]
{
\begin{picture}(0,0)(#1,#2)
\put(-2,-2){$\bullet$}
\setcounter{tempx}{#3}\addtocounter{tempx}{-2}\setcounter{tempy}{#4}\addtocounter{tempy}{-2}
\put(\value{tempx},\value{tempy}){$\bullet$}
\setcounter{tempx}{#5}\addtocounter{tempx}{-2}\setcounter{tempy}{#6}\addtocounter{tempy}{-2}
\put(\value{tempx},\value{tempy}){$\bullet$}
\curvedashes{3,2}\curve(0,0,#3,#4)\curve(0,0,#5,#6)\curve(#3,#4,#5,#6)\curvedashes{}
\setcounter{tempx}{#3/2-5}\setcounter{tempy}{#4/2}\put(\value{tempx},\value{tempy}){\tinyT
#7}
\setcounter{tempx}{(#3+#5)/2+2}\setcounter{tempy}{(#4+#6)/2}\put(\value{tempx},\value{tempy}){\tinyT
#8}
\setcounter{tempx}{#5/2-1}\setcounter{tempy}{#6/2-7}\put(\value{tempx},\value{tempy}){\tinyT
#9}
\end{picture}
}
\newcommand{\squar}[8]
{
\begin{picture}(0,0)(#1,#2)
\put(-2,-2){$\bullet$}
\setcounter{tempx}{#3}\addtocounter{tempx}{-2}\setcounter{tempy}{#4}\addtocounter{tempy}{-2}
\put(\value{tempx},-2){$\bullet$}\put(-2,\value{tempy}){$\bullet$}\put(\value{tempx},\value{tempy}){$\bullet$}

\curvedashes{3,2}\curve(0,0,#3,#4)\curve(0,0,0,#4)\curve(0,0,#3,0)\curve(0,#4,#3,#4)\curve(0,#4,#3,0)\curve(#3,0,#3,#4)\curvedashes{}
\setcounter{tempx}{#3/2}\setcounter{tempy}{#4/2}
\put(\value{tempx},-7){\tinyT #5}\put(-3,\value{tempy}){\tinyT #7}
\put(\value{tempx},\value{tempy}){\tinyT #6}
\setcounter{tempx}{#3/2+2} \put(\value{tempx},\value{tempy}){\tinyT
#8}
\end{picture}
}

\title{On the equi-normalizable deformations of singularities of complex plane curves.}
\author{Dmitry Kerner}
\thanks{Department of Mathematics, Ben Gurion University of the Negev, P.O.B. 653, Be'er Sheva 84105, Israel.
kernerdm@math.bgu.ac.il
\\
The research was constantly supported by the Skirball postdoctoral fellowship of the Center
of Advanced Studies in Mathematics (Mathematics Department of Ben Gurion University, Israel).
\\
Part of the work was done in Mathematische Forschungsinsitute Oberwolfach, during
the author's stay as an OWL-fellow. Some results were published in the preprint \cite{KernerOWLF}}

\begin{document}\setcounter{secnumdepth}{6} \setcounter{tocdepth}{1}\newcounter{tempx}\newcounter{tempy}
\maketitle
\begin{abstract}
We study a specific class of deformations of curve singularities:
the case when the singular point splits to several ones, such that
the total $\de$ invariant is preserved. These are also known as equi-normalizable or equi-generic deformations.
We restrict primarily to the deformations of singularities with
smooth branches.

A natural invariant of the singular type is introduced: the dual
graph. It imposes severe restrictions on the possible
collisions/deformations. And allows to prove some bounds on the
variation of classical invariants in equi-normalizable families.

We consider in details deformations of \omp, the deformations of a singularity into the
collections of \omp s and deformations of the type $x^p+y^{pk}$ into the collections of $A_k$'s.
\end{abstract}
\tableofcontents
\section{Introduction}
In this note we continue the study of collisions/deformations of singular points (started in \cite{Ker07}). For
the relevant notions from singularity theory cf. \cite{AGLV,GLSbook}. Some of the notions are also
 recalled in \S\ref{Sec Known Obstructions}.
The singularity types are named according to Arnol'd's tables \cite[\S I.2.1]{AGLV}.
\subsection{Setup}
Consider plane complex affine singular reduced algebraic curves. Let $\{C_t\}_{t\in T}\subset\mC^2\times T$
be a (flat) family of such curves, with $T$ a small neighborhood (in classical topology) of $0\in\mC^1$.
Let $\si_1..\si_k:T\to\{C_t\}$ be the sections of the family, i.e. the fibre $C_t$ is singular at
the points $\si_1(t)...\si_k(t)$ (for some $t$ the sections might intersect).
We assume that one of the singular points stays at the origin, e.g. $\si_1(t)=(0,0)\in\mC^2$.
Such a family is
called a {\it degeneration} (or degenerating family) if it is not equisingular over $T$ (for local
embedded topological equivalence).
\\
\parbox{14cm}
{\bed A degenerating family is called collision (or split deformation) if the following is satisfied:
\li The family is equisingular over $T\setminus\{0\}$ (for local embedded topological equivalence).
\li The generic fibre  $C_{t\ne0}$ has at least two singular points.
\li The central fibre $C_{t=0}$ has only one singular point.
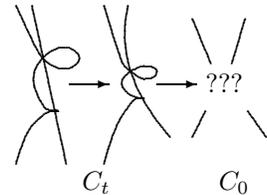
\eed
}
\begin{picture}(0,0)(-20,-5)
\curve(-10,30,-3,15,0,10,4,7,8,6,12,7,14,10,12,13,8,14,4,13,0,10,-3,5,-3,-1,1,-9,4,-10,6,-10)
\curve(6,-10,4,-10,1,-11,-5,-17,-10,-30)\curve(-4,30,9,-30)

\put(15,-40){$C_t$}
\put(10,0){\vector(1,0){15}}
\end{picture}
\begin{picture}(0,0)(-50,-5)
\curve(-1,30,0,5,15,-20)
\curve(-10,30,-3,8,0,5,5,2,10,5,5,7,0,5,-3,2,-3,-1,0,-4,3,-5,5,-5)\curve(5,-5,3,-5,1,-6,-5,-15,-10,-30)
\put(10,0){\vector(1,0){15}}
\end{picture}
\begin{picture}(0,0)(-80,-5)
\curve(10,25,5,8)\curve(11,-10,18,-20)
\curve(-10,25,-3,8)
\curve(-5,-10,-10,-20)

\put(-5,-3){???}
\put(0,-40){$C_0$}
\end{picture}
\\\\
In this paper we always assume the degree of curves to be high enough for the given singularity types
(of the generic and central fibres).

We work with (local embedded topological) singularity types and denote the collisions appropriately.
So $(\mS_1..\mS_k)\to\mS$
means that there exists a flat family of curves $\{C_t\}_{t\in T}$ with the generic fibre $C_{t\ne0}$
having the singularities of the topological types $(\mS_1..\mS_k)$ and the central fibre having only one
singular point of the type $\mS$.

\bex
The simplest cases are $\overset{k}{\underset{i=1}{\coprod}} A_{n_i}\to A_{\sum n_i+k-1}$ and
$\overset{k}{\underset{i=1}{\coprod}} A_{n_i}\cup 2A_1\to D_{\sum n_i+k+2}$  Here the collisions are
possible since the corresponding Dynkin diagrams decompose (see below the general result of \cite{Lyashko83} for ADE).
\eex
The collision $(\mS_1..\mS_k)\to\mS$ is the converse of the deformation $\mS\to(\mS_1..\mS_k)$.
So, the collision is possible if in the versal deformation of some singular germ of type $\mS$ the singularities
of types $(\mS_1..\mS_k)$ appear.

Note that in our case $(\mS_1..\mS_k)\to\mS$ means only that {\it there exists a curve} with the
given types which degenerates to the given type.
In general the possibility of deformation $\mS\to(\mS_1..\mS_k)$ depends not only on the (local embedded)
topological singularity types but on the analytical types (i.e. on the moduli of the types)
 \cite{Pham} \cite{DamonGalligo93} \cite{Jaworski94} and \cite[pg.209]{duPlessis-Wall-book}. So the collision $(\mS_1..\mS_k)\to\mS$
 does not imply that every or generic representative can be degenerated in the prescribed manner.

\bed
A collision/deformation $(\mS_1..\mS_k)\to\mS$ is called prime if it cannot be factorized non-trivially as
$(\mS_1..\mS_k)\to\mS'\to\mS$
(the second arrow being a degeneration of the singularity type).
\eed
The natural questions to ask are:
\\
{\it Given the singularity types $\mS_1..\mS_r$ and $\mS_f$.
Is the collision $(\mS_1..\mS_r)\ra \mS_f$ possible? What are the possible results of prime collisions of $\mS_1..\mS_r$?
What are the possible prime deformations of $\mS_f$?}
\\\\
These questions are highly important in various branches of singularity theory and algebraic geometry, but there
is no hope to get any complete answer at such generality.
\subsection{History}
Some classical (and recent) results are:
\li Every singular point can be deformed to a collection of nodes ($A_1$) in a $\de=const$ way.
(This was claimed in \cite{Albanese28} and (re)proved rigorously in \cite[thm 1.4]{Nobile84}.)
\li every singular point can be deformed to a collection of nodes ($A_1$) and cusps ($A_2$) in
a $\de=const$, $\ka=const$ way \cite[thm 1.1]{DiazHarris88}.
\li Versal deformations of ADE singularities were studied classically.
Let $\mS$ be an ADE singularity. It deforms to a collection $\mS_1..\mS_k$ of (ADE) singularities
iff the collection of Dynkin diagrams of $\mS_1..\mS_k$ can be obtained from that of $\mS$ by
removing some vertices \cite{Lyashko83} (cf. also \cite[\S5.9]{AGLV}). (All the diagrams are taken in the canonical bases of ADE.)
\li Similar result (via Dynkin diagrams in distinguished bases) was obtained in \cite{Jaworski88} (cf. also \cite[\S5.11]{AGLV})
for parabolic singularity types
($P_8=\tilde{E}_6$, $X_9=\tilde{E}_7$, $J_{10}=\tilde{E}_8$).
\li A thorough study of versal deformations of uni-modal types was done in \cite{Brieskorn79},
\cite{Brieskorn81} and \cite{Brieskorn-book} (by studying the properties of Milnor lattices).
\li For many uni-modal and bi-modal types some necessary and sufficient conditions for the decomposition to ADE's
are known. They are formulated in terms of some specific transformations of the canonical Dynkin diagrams.
(cf. \cite{Urabe84} for $J_{3,0},Z_{1,0},W_{1,0}$,$E_{k\le 12}$,$Z_{k\le13},W_{k\le 13}$ and
 \cite{Looijenga81} for $T_{pqr}$).
\li The dependence of some deformations of types $J_{k,0}$ on moduli was studied in \cite[thm 1-thm 4]{Jaworski94} by
checking explicit equations for the singular germs.
\li The deformations of $T_{2pq}$ series (the singularity type of $x^p+x^2y^2+y^q$) has been studied partially.
(Here $X_9=\tilde{E}_7=T_{244}$.) The deformation $T_{244}\to A_7$ was known classically, the
deformation $T_{2pp}\to A_{2p-1}$ was given in \cite[pg.204]{duPlessis-Wall-book}. The deformation
$T_{245}\to A_8$ was constructed in \cite{Stevens04} by brute force computation. The general case
$T_{2pq}\stackrel{?}{\to} A_n$ is open.
\li Some deformations of $E$,$Z$,$Q$ series have been studied in \cite{duPlessis-Wall04}.
\li The question of adjacency of just two singularities (i.e. a singular point deforms to just one singular point)
seems to be more tractable. For the recent advances cf. \cite{Alberich-Roe05}.
\li The case of surfaces in $\mC^3$ is infinitely more complicated (e.g. the whole \cite{Urabe-book}
is a summary of the series of works studying the deformations of just 5 particular singularity types into ADE's).
\\
\\
As it seems, currently no other general results are known. Even worse, it is not clear how to answer (effectively) such
questions in each particular case (except for the case of simple or uni-modal types).
A kind
of brute force computation was given in \cite{Ker07} for a
specific class of singularity types (the so-called linear types). It does not seem to generalize easily to
the arbitrary types.
\\\\
Usually it is very difficult to prove that a deformation exists (e.g. to provide an example).
Rather one seeks for various obstructions. The main classical obstructions are provided by the semi-continuity
of various invariants (cf. \S\ref{Sec Known Obstructions}).
\subsection{The results}
We restrict the consideration to the $\de=const$ deformations, i.e. $\sum\de(\mS_i)=\de(\mS)$,
such that the initial germ (i.e. the central fibre) has all the branches
smooth (\S\ref{Sec delta=const families}).
It is easy to see (proposition \ref{Thm Smooth In branches-->smooth fin branches}) that in such case all
the resulting singulariuties have smooth branches too and $\ka=const$ in the deformation.

The $\de=const$ deformations are particularly important due to their role in the classification
of deformations of sandwiched surface singularities \cite{JongStraten98}.

To any type $\mS$ with smooth branches we associate the {\it dual graph} $\Ga_\mS$.
It is a complete invariant of the local embedded topological type (e.g. is equivalent to the resolution tree).
The importance of the dual graph is due to the classical fact: a $\de=const$ family is equi-normalizable
(theorem \ref{Thm Teissier}). Thus a deformation $\mS\to(\mS_1..\mS_k)$ corresponds to the
decomposition $\Ga_\mS\to\oplus\Ga_{\mS_i}$ (the precise formulation is the theorem \ref{Thm Main Theorem}).
So, the dual graph imposes various restrictions on the possible deformations.
This is our main result.

These restrictions are stronger than some others (e.g. the restriction arising from the Milnor number).
In general they are not weaker than others, e.g they are not implied by the semi-continuity of
the spectrum (cf. example \ref{Example K_5 to 3K_3}).

So we use them all together: those imposed by the dual graph, by the signature of the intersection form on the
middle homology, by the local Bezout theorem, by the spectrum and the Hirzebruch inequality
for \omp (cf. equation (\ref{Eq Hirzebruch Inequality}) below).

As a result we get in many cases necessary conditions not known previously (to the best of author's knowledge).
Below are some consequences of the method (proofs are in \S\ref{Sec Smooth To Smooth Collisions}
\S\ref{Sec Applications}).
\bprop
Let $\mS\to\bigcup\mS_i$ be a $\de=const$ deformation (the types can have singular branches).
Then, the number of branches is bounded: $\bin{r}{2}\le\sum\bin{r_i}{2}$ (the bound is sharp).
\eprop
For types with smooth branches the multiplicity is the number of branches, so this bounds the
change of multiplicity.

The proposition gives an upper bound for the variation of Milnor number, as it satisfies
$\mu-\sum\mu_i=\sum(r_i-1)-(r-1)$.
The lower bound is directly obtained from the classical formulas (\ref{Eq Formulas For Invariants}), and is:
$\frac{\mu}{\sum\mu_i}\le 1+\frac{\sum(r_i-1)-(r-1)}{\sum(p_i-1)^2}$ (the equality is realized for \omp s).

Unfortunately the conditions imposed by the dual graph are not sufficient (cf. remark
\ref{Remark Insufficiency of Dual Graph}).
We do not know whether they can be strengthened in any simple way to become sufficient.
\\
\subsubsection{Results for \omp s}\label{Sec OMP Arrangements Hirzebruch}
 Let $K_p$ denote the topological type of the ordinary multiple point of multiplicity $p$.
 (So $K_2=A_1$, $K_3=D_4$, $K_4=X_9$).
\\{\it Question:} Given the initial type $K_p$ and the final collection $\mS_1..\mS_k$
with $\de(K_p)=\sum\de(\mS_i)$. Is the $\de=const$ deformation $K_{p}\to(\mS_1..\mS_k)$ possible?
It's immediate (proposition \ref{Thm Application to OMP-to-OMP}) that each $\mS_i$ must be an \omp, so the deformation is
$K_{p}\to(K_{p_1}..K_{p_k})$.

Even for such a particular case a satisfactory classification of possible deformations seems to be an open question.
As the deformation is equi-normalizable we can trace each branch separately.
Since all the branches of $K_p$ are transversal and we consider only the germ of curve, can assume they are lines.

If one restricts the question further: assumes the deformed curve consists of lines only, then one has a local arrangement of lines.
So, the classification of the possible
splittings $K_p\to(K_{p_1}...K_{p_k})$ implies in particular the answer to the question:

{\it Given $p$ distinct lines in the plane, which patterns of intersection $(n_2K_2,n_3K_3,...)$ can appear?}

For small $p$ this can be answered by direct classification. In general one can obtain some combinatorial
restrictions (not sufficient).
However there are restrictions of {\it non-combinatorial} origin. In \cite{Hirzebruch83} the Miyaoka-Yau inequality
 for Chern numbers of surfaces was used to prove:
\beq\label{Eq Hirzebruch Inequality}
n_2+\frac{3}{4}n_3\ge p+\sum_{i\ge5} (2i-9)n_i,\text{ provided }n_p=0=n_{p-1}=n_{p-2}
\eeq
We do not know any generalizations or additional non-combinatorial restrictions.
An elementary application of our method leads to the criterion (proved in \S\ref{Sec Applications Deformations of OMP}):
\bprop\label{Intro Thm Application to OMP-to-OMP}
The ($\de=const$) deformation $K_p\to \cupl^k_{i=1} K_{p_i}+\big(\bin{p}{2}-\sum\bin{p_i}{2}\big)K_2$
with $\{p_i\ge k-1\}$ is possible iff $p+\bin{k}{2}\ge\sum p_i$
\eprop
So, for $k\ge5$ this strengthens Hirzebruch's inequality.
Note that for the codimension of the corresponding equisingular strata one has:
\beq
\taues(K_p)-\taues\Big(K_{p_1}..K_{p_k},\big(\bin{p}{2}-\sum\bin{p_i}{2}\big)K_2\Big)=p-\sum p_i+2(k-1)
\eeq
which can be non-positive. In this case the deformation does not exist for the {\it generic} representative of
$K_p$, but only for a very special choice of moduli.
\\\\\\
Another question is: {\it To which collections of \omp s can a given type be deformed in a $\de=const$ way?}
(Generalizing the classical deformation to $\de$ nodes.)

The number of possible scenarios is quite big. One might hope to find some specific prime deformation,
such that all other deformations factorize through this one. Unfortunately this is not the case.
In \S\ref{Sec Applications Canonical Decomp to OMP} we propose a partial result: there always exists
a {\it canonical} deformation into a bunch of \omp s. It is minimal in some sense but it
does not factorize all others.
\\\\
\subsubsection{Decompositions into ADE's}\label{Sec Intro Decomposition into ADE}
An important question is:
To which collections of ADE types deforms a given singularity?

In \S\ref{Sec Applications Deformations into ADE} we study a particular case: $\de=const$ deformations of
the type $x^p+y^{pk}$ (i.e. $p$ smooth branches with equal tangency), denoted by $K_{p,k}$,
 into collections of $A_j$'s. As was noticed above, only $A_j$'s with smooth branches may appear.
We apply the obstructions imposed by the dual graph, the signature of the middle homology
lattice and the spectrum to get:
\bprop
Let $K_{p,k>1}\stackrel{\de=const}{\to}\bigcup n_i A_{2i-1}$. Then $n_{i>k}=0$, $\sum^k_{i=1}in_i=\bin{p}{2}k$ and
$$\sum n_i\ge\frac{(p-1)^2k}{4}+(p-1)-
\Bigg\{\ber\frac{pk}{4}~~p~even\\0~~p~odd,~k~even\\\frac{p-1}{2}~~p~odd,~k~odd\eer$$.

In particular: $n_k\le\frac{(p-1)^2}{4}\frac{k}{k-1}+-\frac{p-1}{k-1}+
\Bigg\{\ber\frac{pk}{4(k-1)}\\0\\\frac{p-1}{2(k-1)}\eer$  and
$2n_1+n_2\ge\frac{(p-1)(p-3)k}{12}+(p-1)-\Bigg\{\ber\frac{pk}{4}\\0\\\frac{p-1}{2}\eer$
\eprop
\subsubsection{Singular branches}
It is not clear how to approach the $\de=const$ deformations when singular branches are present.
An unpleasant fact is: a collection of types with smooth branches can collide ($\de=const$)
to a type with singular branches, such that the collision is prime (i.e. cannot be factorized).
The simplest example is: $A_1+A_3\to D_5$.

When some of the initial branches are singular the dual graph can be still useful. Smoothen the branches
(in a $\de=const$ way): $\mS_i\to\mS^{def}_i$.
Note that in general there are several distinct $\de=const$ branch-smoothings.
So, there are several collections $\{\mS^{def}_i\}$.
Then collide $\{\mS^{def}_i\}$ (in all the possible ways).
This produces some restrictions on the possible results of collision of $\{\mS_i\}\to\mS$ (for $\mS$ with
smooth branches).
\subsection{Acknowledgements}
This work would be impossible without numerous important discussions with E.Shustin and G.M.Greuel.
Many thanks are also to V.Goryunov for important advices.

The main part of this work was done during my stay at the Mathematische Forschungsinsitute Oberwolfach (Germany)
as an OWL-fellow. Many thanks to the staff for excellent working atmosphere.
\section{The classical semi-continuous invariants}\label{Sec Known Obstructions}
Given a singular germ $(C,0)\subset(\mC^2,0)$ and its normalization $\tC\stackrel{\nu}{\to}C$ the following are some simplest topological invariants:
$\mu$ the Milnor number, $mult$ the multiplicity, the $\de:=\dim\nu_*\cO_{\tC}/\cO_C$ invariant
(aka the genus discrepancy, virtual number of nodes etc.), the $\ka$ invariant
(the multiplicity of intersection of the curve with its generic polar), $C=\bigcup C_i$
the branch decomposition, $r$ the number of branches,
$\Ga$ the resolution tree with multiplicities $\{m_i\}$. For the definitions and properties
cf. \cite[I.3.4]{GLSbook}. Some classical formulas are:
\beq\label{Eq Formulas For Invariants}
\de=\frac{\mu+r-1}{2},\hspace{0.5cm}\de=\sum_{p\in \Ga}\frac{m_p(m_p-1)}{2},
\hspace{0.5cm}\de(C)=\sum \de(C_i)+\sum_{i<j}<C_iC_j>,\hspace{0.5cm}\ka=\mu+mult-1
\eeq
For the collection of types $(\mS_1..\mS_k)$ let $\taues(\mS_1..\mS_k)=\suml_i\tau^{es}(\mS_i)$ be
the codimension of the corresponding equsingular stratum in the space of its versal deformation.
Let $Sp(\mS)$ be the spectrum \cite{Steenbrink76} (cf. also \cite[II.8.5]{Kulikov98}).
\bex\label{Example Invariants for K_p and K_pk}
$\star$ Let $K_p$ be the \omp~ of multiplicity $p$. Then:
\beq
\mu(K_p)=(p-1)^2,~~\frac{\ka}{2}=\de(K_p)=\bin{p}{2},~~\taues(K_p)=\bin{p+1}{2},~~
Sp(K_p)=t^{-1+\frac{2}{p}}+2t^{-1+\frac{3}{p}}+3t^{-1+\frac{4}{p}}+...+2t^{1-\frac{3}{p}}+t^{1-\frac{2}{p}}
\eeq
For further use of the semi-continuity we need also the amount of spectral numbers around the origin:
\beq\label{Eq Number Of Spectra Inside Interval For OMP}
\sharp\Big(Sp(K_p)\cap(-\frac{1}{2}+\al,\frac{1}{2}+\al)\Big)=(p-1)^2-\bin{\lfloor(\frac{1}{2}-\al)p\rfloor}{2}
-\bin{\lfloor(\frac{1}{2}+\al)p\rfloor}{2}
\eeq
Here $\bin{n}{2}=0$ for $n<2$.
\\$\star$ Let $\mS$ be the (topological) type of $x^p+y^{q}$ (with $p\le q$). Then:
\beq\ber
\mu(\mS)=(p-1)(q-1),~~\ka(\mS)=(p-1)q,~~\de(\mS)=\frac{(p-1)(q-1)+gcd(p,q)-1}{2},~~
\\\taues(\mS)=\suml^{\lfloor p(1-\frac{2}{q})\rfloor}_{i=2}\lfloor q(1-\frac{i}{p})-1\rfloor,~~
 Sp(\mS)=t^{-1}\suml_{i=0}^{p-2}t^{\frac{i+1}{p}}\suml_{j=0}^{q-2}t^{\frac{j+1}{q}}
\eer\eeq
The amount of the spectral numbers around the origin (now for $q=pk$):
\beq\label{Eq Number Of Spectra Inside Interval For K_pk}
\sharp\Big(Sp(K_{p,k})\cap(-\frac{1}{2}+\al,\frac{1}{2}+\al)\Big)=(p-1)(pk-1)-
\ber\Big(p-\Big\lceil\frac{\lceil pk(\frac{1}{2}+\al)\rceil}{k}\Big\rceil\Big)
\Big(k\frac{p-1+\Big\lceil\frac{\lceil pk(\frac{1}{2}+\al)\rceil}{k}\Big\rceil}{2}-\lceil pk(\frac{1}{2}+\al)\rceil\Big)
\\-
\Big\lfloor\frac{\lfloor pk(\frac{1}{2}+\al)\rfloor}{k}\Big\rfloor
\Big(\lfloor pk(\frac{1}{2}+\al)\rfloor-k\frac{\Big\lfloor\frac{\lfloor pk(\frac{1}{2}+\al)\rfloor}{k}\Big\rfloor+1}{2}
\Big)
\eer\eeq

\eex~
\\\\
Some classical obstructions are:
\li $\mu$, $\ka$, $\de$, $\mu-\de$, $mult$ do not increase under small deformations
(e.g. \cite[theorem 6.1.7]{Buchweitz-Greuel80}). In particular, for $\de=const$ deformation, $(r-1)$ is non-decreasing.
\li If the {\it generic} representative of $\mS$  can be deformed to a curve with $(\mS_1..\mS_k)$
then $\taues(\mS)>\taues(\mS_1..\mS_k)$.
\li The spectrum is semi-continuous in the following sense \cite{Steenbrink85}. Let $Sp(\mS)$ be
the spectrum of the central fibre and $Sp(\cup\mS_i)$ be the joint spectrum of the generic fibre.
Here $Sp(\cup\mS_i)=\sum Sp(\mS_i)$ so that the multiplicities sum up.
For any half-open interval $B_\al=(\al,\al+1]$ let $Sp_B$  be the number of spectral values in the interval
(counting the multiplicities). Then, for every such half-open interval: $Sp_B(\mS)\ge Sp_B(\cup\mS_i)$.
For quasi-homogeneous singularities of curves this holds even for any open interval $(\al,\al+1)$ by \cite{Varchenko83}.

\li The local Bezout theorem can be used as follows. It can be often proved that if a deformation $\mS\to(\mS_1..\mS_k)$
exists then it must be realizable by a curve of small degree. Here one uses the strong criterion of
\cite{Shustin87,GLS96}:
\\
{\it Let $C_d$ be a curve of degree $d$, with the singularity type $\mS$, such that $\tau(\mS)<4d-4$. Then the
  mini-versal deformation of $\mS$ is induced from the parameter space of curves of degree $d$
  (i.e. $|\cO_{\mP^2}(d)|$).}
\\
When this criterion is applicable, one can try to show that a curve (of the small degree) with singularities $\mS_1..\mS_k$
must be reducible and non-reduced, thus forbidding the deformation.

\bex\label{Example X_9} Can an \omp~ of multiplicity 4 (named $X_9$) be deformed to two \omp s of multiplicity 3 (named $D_4$)?
\\
Note that $\taues(X_9)=8=\taues(2D_4)$ so the generic representative of $X_9$ does not deform to $2D_4$.
The simplest classical invariants do not give any restriction: $\mu(X_9)=\tau(X_9)=9>8=\mu(2D_4)=\tau(2D_4)$,
$\de(X_9)=6=2\de(2D_4)$. However, by the previous statement, if the deformation exists it must be
representable by curves of degree 4. But by Bezout theorem a curve of degree 4 with $2D_4$ must have a double
line as a component. So, the deformation $X_9\to2D_4$ is impossible.

The deformation is also prohibited by the semi-continuity of the spectrum
(cf. example \ref{Example Invariants for K_p and K_pk}), Hirzebruch's inequality (\ref{Eq Hirzebruch Inequality})
and the dual graph (proposition \ref{Intro Thm Application to OMP-to-OMP}).
\eex~
\\\\\\
Heavy restrictions arise from the integer cohomology of the Milnor fibre. The cohomology ring is encoded by the lattice
(product in the middle cohomology) \cite[I.1.6]{AGLV}. So the obstruction is:

{\it If a representative of $\mS$ deforms to a curve with singularities $(\mS_1..\mS_k)$ then the direct sum of the cohomology lattices
 of $(\mS_1..\mS_k)$ embeds into the lattice of $\mS$. Correspondingly, there are some bases of vanishing cycles for $(\mS_1..\mS_k)$, $\mS$
 such that the collection of Dynkin diagrams of $(\mS_1..\mS_k)$ is obtained from that of $\mS$ by removing some vertices.}

This restriction is difficult to apply, since it is very difficult to check that one lattice cannot be embedded into another. Alternatively,
one should check that the Dynkin diagram of $\mS$ in all the possible bases cannot be decomposed.
A very painful task even for the ADE types. (Various results on the behavior of Dynkin diagrams under the base change
and specific criteria can be found e.g. in \cite{Brieskorn-book,Urabe-book}.)

In this work we use this obstruction only partially by the signature of the middle homology. Namely, consider the stabilization of
the curve singularities to surfaces, i.e. instead of the curve $f(x,y)=0$ one has the suspension $f(x,y)+z^2=0$, the surface singularity in $\mC^3$.
Then the intersection form in the middle (co)homology is symmetric. And if the signature forbids a collision for surfaces then so is
for curves.

Given a lattice $(L,<,>)$ form the
corresponding vector space $V=L\underset{\mZ}{\otimes}\mR$ with the induced quadratic form $Q:~V\to\mR$ (by $v\to<v,v>$). Thus if $L_1\subset L_2$
one has $V_1\subset V_2$ and $Q_1=Q_2|_{V_1}$. This brings a restriction on the signature:
\bprop $\bullet$
Let $(\mu_+(Q_i),\mu_0(Q_i),\mu_-(Q_i))$ be the number of positive/zero/negative eigenvalues. Then $\mu_+(Q_1)\le\mu_+(Q_2)$,
$\mu_-(Q_1)\le\mu_-(Q_2)$ and $\mu_0(Q_1)\le\mu_0(Q_2)+dim(V_2)-dim(V_1)$.
\li Let $Q^{(N)}_i=Q_i\mod N$ be the quadratic form reduced modulo N. Then:
 $\mu_0(Q^{(N)}_1)\le\mu_0(Q^{(N)}_2)+dim(V_2)-dim(V_1)$.
\eprop
\bpr
Here only the bound on $\mu_0$ possibly needs an explanation. Fix any base for $L_1$ and extend it to a base of $L_2$.
Let $Q_1,Q_2$ be  the matrices of the intersection forms in this basis. So, $Q_1$ is a submatrix of $Q_2$: $Q_2=\bpm Q_1&A\\B&C\epm$.
Apply the conjugation $Q_2\to U_1Q_2U_2$ with $U_i\in GL(V_2)$ preserving $V_1\subset V_2$ (so that the block structure is preserved).
By such conjugation can "diagonalize" the blocks $A,B$:
\beq
Q_2=\bpm Q_1&\bpm *&0\\0&*\\0&0\epm\\\bpm *&0&0\\0&*&0\epm&**  \epm
\eeq
Hence $\dim Ker(Q_1)-\dim Ker(Q_1)\le dim(V_2)-dim(V_1)$
\epr

The numbers $\mu_\pm,\mu_0$ can be computed from the spectrum: $\mu_-=\sharp(Sp\cap(-\frac{1}{2},\frac{1}{2})$
and $\mu_+=2\times\sharp(Sp\cap(-1,-\frac{1}{2}))$.

For the quasi-homogeneous case they can be calculated as follows \cite{Steenbrink77}.

Let $\{f(x_1,x_2,x_3)=0\}\subset\mC^3$ be a quasi-homogeneous surface singularity, with $x_i$ of weight $w_i$
and $f$ of weight 1. Take the monomial basis for its  Milnor algebra $\mC[x_1,x_2,x_3]/j_f$: $e_1...e_\mu$.
For each such monomial define the weight function $l(x^{n_1}_1..x^{n_3}_3):=\sum w_i(n_i+1)$.
\bprop\cite[theorem 2]{Steenbrink77}Given the germ of a quasi-homogeneous surface $\{f=0\}\subset(\mC^3,0)$.
Let $M_+\oplus M_0\oplus M_-$ be the decomposition of the (co)homology lattice according to the signature of
the intersection product. Then the spaces are spanned by the residue forms:
\beq
M_0=\{e_i\frac{dxdydz}{df}|l(e_i)\in\mZ\},~~
M_+=\{e_i\frac{dxdydz}{df}|l(e_i)\notin\mZ,~\lfloor l(e_i)\rfloor\text{ is even}\},~~
M_-=\{e_i\frac{dxdydz}{df}|l(e_i)\notin\mZ,~\lfloor l(e_i)\rfloor\text{ is odd}\}
\eeq
In particular $\mu_+=|M_+|$, $\mu_0=|M_0|$, $\mu_-=|M_-|$.
\eprop
\bex
Let $\mS$ be the topological type of $x^p+y^{pk}+z^2$ (for $k\ge1$). Direct calculation gives:
\li $p$ even: $\mu_0=p-2$, $\mu_+=\frac{p-2}{2}(\frac{pk}{2}-2)$ and  $\mu_-=\frac{3p-2}{4}kp-(p-1)$
\li $p$ odd, $k$ even: $\mu_0=p-1$, $\mu_+=(\frac{p-1}{2})^2k-(p-1)$ and  $\mu_-=\frac{(p-1)(3pk+k-4)}{4}$
\li $p$ odd, $k$ odd: $\mu_0=0$, $\mu_+=(\frac{p-1}{2})^2k-\frac{p-1}{2}$ and  $\mu_-=\frac{(p-1)(3pk+k-2)}{4}$
\eex
\bex\label{Example J_10} Can the singularity of type $J_{10}$ (locally $x^3+\al xy^4+y^6$, with $\al$ the modulus)
be deformed to three tacnodes ($3A_3$)?
Since $\taues(J_{10})=9=\taues(3A_3)$ the generic representative cannot be deformed.
Other invariants are: $\mu(J_{10})=\tau(J_{10})=10>\mu(3A_4)=\tau(3A_3)$, $\de(J_{10})=6=\de(3A_3)$, $\ka(J_{10})=12=\ka(3A_3)$.

The local Bezout theorem does not give any restriction since the minimal degree of curve possessing $J_{10}$ is 5.

The deformation is forbidden by the signature of the middle homology. In fact for $J_{10}$ have
$(\mu_+\mu_0,\mu_-)=(0,2,8)$, for $3A_3$ have
$(\mu_+\mu_0,\mu_-)=(0,0,9)$.
(Recall that the signatures are calculated for the stabilizations: surfaces in $\mC^3$.)
\eex
\bex\label{Example K_5 to 3K_3}$\bullet$ Does the deformation $K_5\to3 K_3+A_1$ exist? (Can a representative
of the \omp~of multiplicity $5$ be deformed to 3 points of types $D_4$ and one node?)
The invariants as above do not forbid it:
\beq\ber
\de(K_5)=10=3\de(K_3)+\de(A_1),~~\mu_+(K_5)=2,~~\mu_0(K_5)=0,~~\mu_-(K_5)=14>13=3\mu(K_3)+\mu(A_1),\\
Sp(K_5)=t^{-\frac{3}{5}}+2t^{-\frac{2}{5}}+3t^{-\frac{1}{5}}+4t^0+3t^{\frac{1}{5}}+2t^{\frac{2}{5}}+t^{\frac{3}{5}},~~~~
Sp(3K_3\cup A_1)=3t^{-\frac{1}{3}}+7t^0+3t^{\frac{1}{3}}
\eer\eeq
Hirzebruch's inequality (\ref{Eq Hirzebruch Inequality}) is irrelevant (since $n_{p-2}\neq0$).
\li Let $K_{p,k}$ be the singularity of the type $x^p+y^{pk}$. Does the deformation $K_{4,3}\to 7A_3+4A_1$ exist?
It is not forbidden by the classical invariants above.
\eex
However, in both cases the deformations are forbidden by the dual graph
(cf. propositions \ref{Thm Application to OMP-to-OMP} and \ref{Thm K_pk to A_n's}).

\section{The $\de=const$ families.}\label{Sec delta=const families}
In this note rather than restrict the singularity types, we restrict to the $\de=const$ deformations.
As the deformations are equi-normalizable the problem consists of two parts:
\ls To understand the $\de=const$ deformations of branches
\li To understand the ways to combine possible deformations of the branches to a deformation of the curve.
\\
Note, this doesn't imply that any deformation can be factorized as
$\cup(C_i,0)\to\cup(C^{def}_i,0)\to\cup(C^{def}_i,0_i)$,
i.e. first each branch is deformed while preserving the singularity at the origin, then the branches
are moved (but their types are preserved). The simplest example is $A_4\to A_2+A_1$.

We consider the second question. In fact we restrict  further:
\bprop\label{Thm Smooth In branches-->smooth fin branches}
Suppose the deformation $\mS\stackrel{\de=const}{\to}(\mS_1..\mS_k)$ is possible.
 If all the branches of $\mS$ are smooth then all the branches of $(\mS_1..\mS_k)$ are
smooth and the deformation is $\ka=const$.
If at least one branch of $\mS$ is non-smooth and the deformation is $\ka=const$,
then at least one branch of $\cup\mS_i$ is non-smooth.
\eprop
\bpr  The $\de=const$ means, cf. eq. (\ref{Eq Formulas For Invariants}): $\mu+r-1=\sum(\mu_i+r_i-1)$.
The semi-continuity of $\ka$ means: $\mu+p-1\ge\sum(\mu_i+p_i-1)$ (for $p_i=mult(\mS_i)$).
Thus $0=p-r\ge\sum(p_i-r)$ implies $p_i=r_i$ i.e. the branches are smooth.
\epr
\\
The $\de=const$ deformations/degenerations are called equi-generic, the $\de=const$, $\ka=const$
are called equi-classical \cite{DiazHarris88}.

The $\de=const$ deformations/degenerations are very special: they are equinormalizable.
\bthe\label{Thm Teissier}\cite{Teissier}, cf. also \cite[Theorem I.2.54]{GLSbook}.
Let $S=\{C_t\}$ be a family of plane curves (considered as a fibred surface). Let $\tS\ra S$ be
the normalization of the surface and $\tC_0\ra C_0$ the corresponding map of the central fibres.
Then $\de(C_0)-\de(C_{t\ne0})=\de(\tC_0)$. In particular the family is equinormalizable (i.e. $\tC_0$ is smooth)
iff it is $\de=const$.
\ethe
This key property enables to associate a new obstruction: the dual graph.
%
%
%
%
%
%
\section{The dual graph}\label{Sec Dual Graph}
\bex\label{Example Collision of Two Nodes To Tacnode}
As a motivation consider the ($\de=const$) collision $A_{2k-1}+A_{2l-1}\to A_{2k+2l-1}$.
\\
\parbox{13cm}
{Let $S=\{C_t\}_t\to T$ be the fibred surface formed by the degenerating family. $S$ has non-isolated singularities.
Consider the normalized surface $\tS=\{\tC_t\}_t\to T$.
Since the collision is $\de=const$, the surface is smooth and each fibre is the normalization of the corresponding curve.
Let $p_1,p_2$ and $q_1,q_2$ be the pre-images of $A_{2k-1},A_{2l-1}$ (in the generic fibre).
As $p_1,p_2$ are glued by the normalization connect them by a dotted line (with multiplicity $k$,
the local intersection of the branches). Similarly for $q_1,q_2$. Then the collision can
be traced on the normalized surface as the addition of graphs: $q_i\ra p_i$, the edges merge
(and their weights are added). Pictorially:
$\bullet\stackrel{k}{-}\bullet+\bullet\stackrel{l}{-}\bullet=\bullet\stackrel{k+l}{--}\bullet$
}
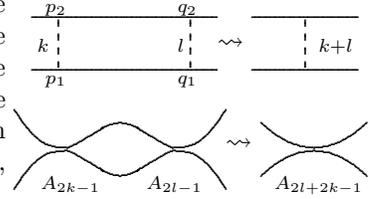
\begin{picture}(0,0)(0,20)
\curve(0,15,10,3,20,0.5,25,3,33,8,40,10,47,8,55,3,60,0.5,70,3,80,15)
\curve(0,-15,10,-3,20,-0.5,25,-3,33,-8,40,-10,47,-8,55,-3,60,-0.5,70,-3,80,-15)\put(80,0){$\rightsquigarrow$}
\put(10,-15){$\tinyA A_{2k-1}$}  \put(50,-15){$\tinyA A_{2l-1}$}
\end{picture}
\begin{picture}(0,0)(-90,20)
\curve(0,-10,20,-0.5,40,-10)\curve(0,10,20,0.5,40,10)  \put(5,-15){$\tinyA A_{2l+2k-1}$}
\end{picture}
\begin{picture}(0,0)(0,-10)
\curve(0,20,70,20) \curve(0,0,70,0)\curvedashes{3,2}\curve(10,22,10,0) \curve(60,22,60,0)\curvedashes{}
\put(5,-5){$\tinyA p_1$} \put(5,22){$\tinyA p_2$}\put(2,7){$\tinyA k$}   \put(55,22){$\tinyA q_2$} \put(55,-5){$\tinyA q_1$} \put(55,7){$\tinyA l$}
\put(70,8){$\rightsquigarrow$}
\end{picture}
\begin{picture}(0,0)(-80,-10)
\curve(0,20,40,20) \curve(0,0,40,0)\curvedashes{3,2}\curve(20,22,20,0)\curvedashes{}\put(25,7){$\tinyA k+l$}
\end{picture}
\eex

\bed
For the singularity $\mS$ with smooth branches the dual graph $\Ga_\mS$ is a complete graph with weighted edges, such that:
\li The vertices $\{v_i\}$ of $\Ga_\mS$ correspond bijectively to the branches $\{C_i\}$ of $\mS$
\li Let the vertices $v_i,v_j\in\Gamma_\mS$ correspond to the branches $C_iC_j$. The weight of the edge $\overline{v_iv_j}$ equals $<C_i,C_j>$.
\\ If $C_i\pitchfork C_j$ (i.e. $<C_iC_j>=1$) then the weight is omitted.
\eed
\bex Below are some graphs for the \omp~of multiplicity $p$ (denoted by $K_p$), for the type of $x^p+y^{pk}$ (denoted by $K_p(k)$) etc.
\\
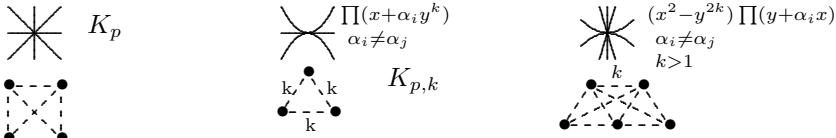
\begin{picture}(0,0)(-20,10)
\curve(0,10,0,-10)\curve(10,10,-10,-10) \curve(-10,10,10,-10) \curve(10,0,-10,0)
\squar{10}{40}{20}{20}{}{}{}{}
\put(20,0){$K_p$}
\end{picture}
\begin{picture}(0,0)(-110,10)
\curve(0,10,10,0,20,10)\curve(0,-10,10,0,20,-10)\curve(0,0,20,0)
\put(21,5){$\tinyA \prod(x+\al_i y^k)$}  \put(25,-4){$\tinyA \al_i\ne\al_j$}
\triang{0}{30}{10}{15}{20}{0}{k}{k}{k}
\put(40,-20){$K_{p,k}$}
\end{picture}
\begin{picture}(0,0)(-220,10)
\curve(0,5,10,0,20,5)\curve(0,-5,10,0,20,-5)\curve(10,10,10,-10)
 \curve(7,10,13,-10) \curve(7,-10,13,10)
\put(25,5){$\tinyA (x^2-y^{2k})\prod(y+\al_i x)$}  \put(28,-4){$\tinyA \al_i\ne\al_j$}\put(28,-13){$\tinyA k>1$}
\end{picture}
\begin{picture}(0,0)(-220,30)
\put(-2,-2){$\bullet$} \put(18,-2){$\bullet$}
\put(-12,-17){$\bullet$} \put(8,-17){$\bullet$} \put(28,-17){$\bullet$}
 \curvedashes{3,2}\curve(0,0,20,0)\curve(-10,-15,30,-15)
 \curve(0,0,-10,-15) \curve(0,0,10,-15) \curve(0,0,30,-15)
 \curve(20,0,-10,-15) \curve(20,0,10,-15) \curve(20,0,30,-15)\curvedashes{}
\put(8,3){$\tinyA k$}
\end{picture}
\eex
~\\\\\\\\
Though defined for a germ, the dual graph depends on the (local embedded topological) singularity type only.
\bprop\label{Thm Delta of a Graph}Let $(C,0)$ be a germ with smooth branches.
The dual graph $\Ga_{\mS(C)}$ is well defined (independent of the representative of $\mS$) and is a complete invariant of the
local embedded topological type $\mS(C,0)$. In particular:
\li the number of branches=the multiplicity=the number of vertices of the graph
\li $\frac{\ka(\mS)}{2}=\de(\mS)=\sum w(\overline{v_iv_j})$ (the sum is over all the edges).
\eprop
\bpr
$\bullet$$\Ga_\mS$ is well defined because the number of branches and the multiplicities of their pairwise intersections are topological invariants.
Conversely, from this data the resolution tree is immediately restored.
\li Use $\de(\mS)=\sum_{i<j}<C_i,C_j>$, cf. (\ref{Eq Formulas For Invariants})
\epr
The weights of the dual graph satisfy some consistency conditions:
\bel\label{Thm Necessary Conditions for a dual graph}
For any path of vertices $<v_1...v_n>$ in the graph one has:
$w(\ol{v_1v_n})\ge\min\Big(w(\ol{v_1v_2})w(\ol{v_2v_3})..w(\ol{v_{n-1}v_n})\Big)$.
In particular, let $w=w(\ol{v_iv_j})$ be the minimal among the weights of edges in the graph. Then for any vertex
$v_n$ either $w(\ol{v_iv_n})=w(\ol{v_iv_j})$ or $w(\ol{v_jv_n})=w(\ol{v_iv_j})$.
Moreover, let $\Ga'$ be the graph obtained from the dual graph by removing all the edges of weight $\le k$.
Then $\Ga'$ is the (disjoint) union of complete graphs (possibly isolated points).
\eel
This follows from the observation $\min(<C_iC_j>,<C_jC_n>)\le <C_iC_n>$
(immediate from the resolution tree) and consideration of all the paths $<v_iv_jv_n>$.

The converse is true: this condition is sufficient for a graph to be the dual graph.
\bel\label{Thm Sufficient Conditions for a dual graph}
$\bullet$ Let $\Ga$ be a complete graph with weighted edges, such that for each triple of vertices $v_iv_jv_k\in\Ga$
the weights satisfy (possibly after a permutation): $w(\ol{v_iv_j})=w(\ol{v_iv_k})\le w(\ol{v_kv_j})$. Then $\Ga$
is the dual graph for some topological type, i.e. $\Ga=\Ga_\mS$.
\li Let $\Ga$ be the dual graph of some singularity. Let $\Ga_1\subset\Ga$ be any full complete subgraph (i.e.
if $v_i,v_j\in\Ga_1\subset\Ga$ then $w(\ol{v_iv_j})_{\Ga_1}=w(\ol{v_iv_j})_\Ga$). Then $\Ga_1$ is also
the dual graph of some singularity.
\eel
\bpr
We construct an explicit representative of $\mS$. To each vertex $v_i$ associate an abstract smooth branch-germ
$(C_i,0)\approx(\mC^1,0)$ and embed them into $\mC^2$ inductively. Start from $i_1:(C_1,0)\hookrightarrow(\mC^2,0)$.
Let $i_2:(C_2,0)\hookrightarrow(\mC^2,0)$ be an embedding such that $<i_1(C_1)i_2(C_2)>=w(\ol{v_1v_2})$,
and is generic otherwise. Suppose the branches $C_1..C_k$ are embedded, such that $<i_j(C_j)i_l(C_l)>=w(\ol{v_jv_l})$
(for $1\le j<l\le k$).

For $v_{k+1}$ and $C_{k+1}$ consider the integers $\{w(\ol{v_{k+1}v_j})\}_{j=1..k}$.
Suppose the maximum is obtained for $w(\ol{v_{k+1}v_l})$ (if such $l$ is non-unique, choose any of them).
Embed $i_{k+1}:(C_{k+1},0)\hookrightarrow(\mC^2,0)$
such that $<i_{k+1}(C_{k+1})i_l(C_l)>=w(\ol{v_{k+1}v_l})$ but is generic otherwise.

Then for any $j\le k$:  $<i_{k+1}(C_{k+1})i_j(C_j)>\le<i_{k+1}(C_{k+1})i_l(C_l)>$. So, by the assumption
of the proposition: $<i_{k+1}(C_{k+1})i_j(C_j)>=<i_l(C_l)i_j(C_j)>\le<i_{k+1}(C_{k+1})i_l(C_l)>$.
But $<i_l(C_l)i_j(C_j)>=w(\ol{v_lv_j})\le w(\ol{v_lv_{k+1}})$, thus $<i_{k+1}(C_{k+1})i_j(C_j)>=w(\ol{v_jv_{k+1}})$.
And this proves the criterion.

The second statement is obvious.
\epr

The dual graphs can be often added and decomposed.
\bed\label{Def Addition Of The Dual Graphs}
The (complete weighted) graph $\Ga_\mS$ decomposes into the union $\sqcup_i \Ga_{\mS_i}$ (write $\Ga_\mS=\oplus\Ga_{\mS_i}$) if there
is a map $\phi:\sqcup_i \Ga_{\mS_i}\to\Ga_\mS$ (vertices-to-vertices, edges-to-edges) surjective on vertices and edges of $\Ga_\mS$  such that:
\li For any two distinct vertices $v_1,v_2\in\Gamma_{\mS_i}$: $\phi(v_1)\ne\phi(v_2)\in\Ga_\mS$ and $\overline{\phi(v_1),\phi(v_2)}=\phi(\overline{v_1,v_2})$
\li The weights of edges add up, i.e. for $v_1,v_2\in\Ga_\mS$: $w(\overline{v_1,v_2})=\underset{l\in\phi^{-1}(\overline{v_1,v_2})}{\sum} w(l)$
\eed
\bex $\bullet$ Every dual graph can be decomposed to the union of $K_2$'s (two vertices and an edge of weight 1).
This corresponds to the standard deformation of the the singularity to $\de$ nodes.
\li Let $K_p$ be the complete graph on $p$ vertices (all weights are one).
 Direct check shows that $K_4$ cannot be decomposed to the union of two $K_3$'s.
 This corresponds to the impossibility of deformation $X_9\to2D_4$.
\eex
Sometimes the dual graphs can be also subtracted.
\bed\label{Def Dual Graph Subtraction}
Suppose there exists an embedding $\Ga_1\into\Ga_2$ of vertices and edges with non-decreasing weights,
i.e. for any two vertices $v_kv_j\in\Ga_1$ have: $w(\ol{v_kv_j})\le w(\ol{i(v_k)i(v_j)})$.
Denote by $\Ga_2-i(\Ga_1)$  the weighted graph obtained from $\Ga_2$ by subtracting
the weights: $w(\ol{i(v_k)i(v_j)})-w(\ol{v_kv_j})$ (for edges of $\Ga_2$ in the image of $\Ga_1$).
If $w(\ol{i(v_k)i(v_j)})=w(\ol{v_kv_j})$ then the edge  $\ol{i(v_k)i(v_j)}$ is erased.
All the isolated vertices are erased too.
\eed
In general $\Ga_2-i(\Ga_1)$ can have several connected components.\\
\bex
\triang{-10}{0}{10}{15}{20}{0}{k}{k}{2k}\hspace{1.3cm}-\triang{-10}{0}{10}{15}{20}{0}{k}{k}{k}
\hspace{1.2cm}= $\bullet\stackrel{k}{\cdot\cdot}\bullet$
This can be written also as:\\\\
\beq\label{Eq Addition Subtraction
Dual Graph}
\triang{-10}{0}{10}{15}{20}{0}{k}{k}{2k}\hspace{1.3cm}=\triang{-10}{0}{10}{15}{20}{0}{k}{k}{k}
\hspace{1.2cm}+\bullet\stackrel{k}{\cdot\cdot}\bullet=
\oplusl_{i=1}^k\triang{-10}{0}{10}{15}{20}{0}{}{}{}\hspace{1.2cm}
+\oplusl^k_{i=1}\bullet\stackrel{}{\cdot\cdot}\bullet
\eeq \eex
The last example is a particular case of an important operation: {\it the canonical
decomposition} into \omp s.

Let $w$ be the minimal among the weights of edges of $\Ga_\mS$. Let $K_p$ be the complete graph on $|\Ga_\mS|$
vertices, with each edge of weight $1$. Consider the graph $\Ga_\mS-w K_p$.
 It is immediate from the proposition \ref{Thm Sufficient Conditions for a dual graph}
that each connected component of $\Ga_\mS-w K_p$ is the dual graph of some singularity type, i.e.
$\Ga_\mS-w K_p=\oplus\Ga_{\mS_i}$. Apply the same procedure to each of $\Ga_{\mS_i}$, till one gets the
collection of graphs with edges of weight 1.
\bed\label{Def Canonical decomposition of dual graph}
The so defined decomposition $\Ga_\mS=\oplus w_i K_{p_i}$ is called
the canonical decomposition for the type $\mS$.
\eed
An example of canonical decomposition is the last equality in equation (\ref{Eq Addition Subtraction Dual Graph}).
Note that the canonical decomposition is well defined and depends on the initial singularity type only.

\beR\label{Rem Dual Graph Sing Branch}
 It is not clear how to define the dual graphs for types with singular branches. A trivial choice is,
to ignore the singularities of branches and trace the contacts of branches only (i.e. the minimal number of
blowups needed to separate them). Eventhough one looses lots of information, this leads to a new
semi-continuous invariant (cf. corollary \ref{Thm Bound On r for any type}).
\eeR
\section{The dual graph restrictions on $\de=const$ deformations}\label{Sec Smooth To Smooth Collisions}
Consider $\de=const$ family of plane curves: $\{C(t)\}_{t\in I}$. Let the central fibre possess only one singular
point (at $0\in\mC^2$), with smooth branches. Let the generic fiber possess the
singularities at $\{x_i\in\mC^2\}$ (by the proposition \ref{Thm Smooth In branches-->smooth fin branches}
all the branches are smooth).

Denote the result of the deformation as the map of germs $\bigsqcup(C_i,x_i)\ra(C,0)$.
Let $\{\Ga_{\mS_i}:=\Ga(C_i,x_i)\}$  and $\Ga_\mS:=\Ga(C,0)$ be the dual graphs.
\bthe\label{Thm Main Theorem}
The $\de=const$ deformation $(C,0)\to\bigsqcup(C_i,x_i)$  (all the branches are smooth)
induces the decomposition $\Ga_\mS=\oplus\Ga_{\mS_i}$.
\ethe
\bpr We should construct the needed map (cf. definition \ref{Def Addition Of The Dual Graphs}).
Use Teissier's theorem \ref{Thm Teissier}.
\\
{\bf Construction of the map}.
Let $\{\tC(t)\}_{t\in I}\ra\{C(t)\}_{t\in I}$ be the normalization of the family.
Namely, for every $t$, the curve  $\tC(t)$ is smooth and the map $\tC(t)\stackrel{\nu_t}{\ra} C(t)$
is the normalization. Let $(C,0)=\bigcup(C^\al,0)$ be the branch decomposition and
$\bigsqcup(\tC^\al,0^\al)\stackrel{\bigcup \nu_\al}{\to}\bigcup(C^\al,0)$ the corresponding preimages.

As the singular points collide, the preimages of each of them in the normalized family $\{\tC_t\}_{t\in I}$ converge
 to one of the points $0^\al$.
This defines the map $\pi_i$ of the vertices of each graph $\Ga(C_i,x_i)$ to the graph of the central fibre $\Ga(C,0)$.

{\bf Properties of the map}.
Let $(C^\al,0)$ be a (smooth) branch of the central fibre, let $(\tC^\al,0^\al)$ be the corresponding branch
on the normalized surface. Let $\bigsqcup \tC^\al_i(t)$ be the normalization of all the branches that
merge to $\tC^\al$. Correspondingly the (smooth) branches $\bigcup C^\al_i$ specialize to $C^\al$.
But $C^\al$ is smooth (and reduced), thus the branches $\{\nu(\tC^\al_i)\}_i$
cannot intersect.

For the corresponding graph this means: let $v^\al_i,v^\al_j$ be two vertices of $\bigsqcup_k\Ga_k$,
sent by $\bigcup \pi_k$ to the same vertex $v^\al\in\Ga_\mS$. Then they are not connected (i.e. belong to distinct
complete subgraphs).  Therefore no edge of $\bigsqcup_k\Ga_k$ is contracted by $\bigcup \pi_k$.

Similarly, note that $<C^\al C^\be>\ge<\Big(\bigsqcup\nu(\tC^\al_i(t))\Big)\Big(\bigsqcup\nu(\tC^\be_j(t))\Big)>$.
But for the singular point with smooth branches $\de(\bigcup C^\al)=\sum<C^\al  C^\be>$ (with the sum over all
pairs of branches). Correspondingly (for $v^\al,v^\be\in\Ga_\mS$):
\beq
w(\ol{v^\al v^\be})=<C^\al C^\be>=<\Big(\bigsqcup\nu(\tC^\al_i(t))\Big)\Big(\bigsqcup\nu(\tC^\be_j(t))\Big)>=
\sum_{\ber\pi(v^\al_i)=v^\al\\\pi(v^\be_j)=v^\be\eer} w(\ol{v^\al_iv^\be_j})
\eeq
Finally the map is surjective as no new branches are created in $\de=const$ collisions.
\epr
\bex
To illustrate the use of the last proposition, consider the collision of $\de$ nodes. So we have $\de$ graphs
(each being just an edge with two vertices). From these building blocks we should glue a complete
graph with weighted edges (such that the weights are added). Below are some examples for low $\de$.
\beq\ber
\oplusl^2_{i=1}\bpm \bullet\\|\\\bullet\epm\rightsquigarrow \bullet\overset{2}{\underset{A_3}{\cdot\cdot\cdot}}\bullet
\hspace{0.5cm}
\oplusl^3_{i=1}\bpm \bullet\\|\\\bullet\epm\rightsquigarrow \bullet\overset{3}{\underset{A_5}{\cdot\cdot\cdot}}\bullet~~\rm{or}~~~~~
\ber\triang{0}{0}{10}{10}{20}{0}{}{}{}\\~~D_4\eer
\hspace{0.5cm}
\oplusl^4_{i=1}\bpm \bullet\\|\\\bullet\epm\rightsquigarrow \bullet\overset{4}{\underset{A_7}{\cdot\cdot\cdot}}\bullet~~\rm{or}~~~
\ber \triang{0}{0}{10}{10}{20}{0}{}{}{2}\\\\~~D_6\eer
\hspace{0.5cm}
\oplusl^5_{i=1}\bpm \bullet\\|\\\bullet\epm\rightsquigarrow \bullet\overset{5}{\underset{A_9}{\cdot\cdot\cdot}}\bullet~~\rm{or}~~~
\ber\triang{0}{0}{10}{10}{20}{0}{}{}{3}\\\\~~D_8\eer
\\\\
\oplusl^6_{i=1}\bpm \bullet\\|\\\bullet\epm\rightsquigarrow \bullet\overset{6}{\underset{A_{11}}{\cdot\cdot\cdot}}\bullet~~\rm{or}~~~
\ber\triang{0}{0}{10}{10}{20}{0}{}{}{4}\\\\~~D_{10}\eer\hspace{0.15cm}\rm{or}~~~
\ber\triang{0}{0}{10}{10}{20}{0}{2}{2}{2}\\\\~~J_{10}\eer\hspace{0.15cm}\rm{or}~~
\ber\squar{0}{0}{15}{15}{}{}{}{}
\\X_9\eer
\hspace{2cm}
\oplusl^7_{i=1}\bpm \bullet\\|\\\bullet\epm\rightsquigarrow \bullet\overset{7}{\underset{A_{13}}{\cdot\cdot\cdot}}\bullet~~\rm{or}~~~
\ber\triang{0}{0}{10}{10}{20}{0}{}{}{5}\\\\~~D_{12}\eer\hspace{0.15cm}\rm{or}~~~
\ber\triang{0}{0}{10}{10}{20}{0}{2}{2}{3}\\\\~~J_{2,2}\eer\hspace{0.15cm}
\rm{or}~~\ber\squar{0}{0}{15}{15}{2}{}{}{}
\\\\~~X_{1,2}\eer
\eer\eeq
Note that in the case of 5 nodes we discard the possibility \triang{0}{0}{10}{10}{20}{0}{2}{2}{}
\hspace{0.8cm} as this graph is not a dual graph of any singularity
(cf. proposition \ref{Thm Necessary Conditions for a dual graph}).
Similarly for 6 nodes we discard the case \triang{0}{0}{10}{10}{20}{0}{1}{2}{3}\hspace{0.8cm}, while for 7 nodes
we discard \triang{0}{0}{10}{10}{20}{0}{1}{3}{3}\hspace{0.8cm} and \triang{0}{0}{10}{10}{20}{0}{1}{2}{4}\hspace{0.8cm}
\eex~
\\\\
The obstruction imposed by the dual graph is stronger than some others and in particular provides a bound on the jump of the
Milnor number.
\bprop
Let $(\mS_1..\mS_k),\mS$ be the types with smooth branches such that the dual graph decomposes: $\Ga_\mS=\oplus \Ga_{\mS_i}$.
Then $\de(\mS)=\sum\de(\mS_i)$,  $\ka(\mS)=\sum\ka(\mS_i)$,  $mult(\mS)\ge\max\Big( mult(\mS_i)\Big)$, $\mu_\mS>\sum_i\mu_{\mS_i}$.
Moreover $\mu_\mS-\sum_i\mu_{\mS_i}=\sum(r_i-1)-(r-1)$, where $\sum\bin{r_i}{2}\ge\bin{r}{2}$.
\eprop
\bpr
Here only the bound on the Milnor number should be explained. The equality $\mu_\mS-\sum_i\mu_{\mS_i}=\sum(r_i-1)-(r-1)$ arises
from $\de(\mS)=\sum\de(\mS_i)$. The bound $\sum\bin{r_i}{2}\ge\bin{r}{2}$ arises from the surjectivity of the map $\sqcup\Ga_{\mS_i}\to\Ga_\mS$.
\epr
It seems that the obstruction imposed by the dual graph is not implied by any known obstructions,
in particular it is not weaker than the spectrum (cf. example \ref{Example K_5 to 3K_3}).
\beR\label{Remark Insufficiency of Dual Graph}
Unfortunately the necessary condition from the dual graph is far from being sufficient.
For example, the deformation $J_{10}\to3A_3$ is impossible (cf. example \ref{Example J_10}), but the
corresponding graph certainly decomposes:\\
\beq
\triang{-10}{0}{10}{15}{20}{0}{3}{3}{3}\hspace{1.3cm}\to\oplusl^3\triang{-10}{0}{10}{15}{20}{0}{}{}{}
\eeq
The conditions imposed by the dual graph are not sufficient even for deformations of an \omp~ into \omp s.
Indeed, there are classical examples of the decompositions of complete graphs: $K_p\to\oplus K_{p_i}$
with each $p_i>3$. Contrary to Hirzebruch's inequality (\ref{Eq Hirzebruch Inequality}) for arrangement of lines.
\eeR
\subsection{The case of non-smooth branches in the initial types.}\label{Sec Non Smooth Initial Branches}
For the case of non-smooth branches the $\de=const$ deformation is still equi-normalizable
(theorem \ref{Thm Teissier}). So, one can consider the dual graph (cf. \ref{Rem Dual Graph Sing Branch})
and the idea of the proof of the proposition \ref{Thm Main Theorem}
gives a weaker statement:
\bprop
Any $\de=const$ deformation $\mS\to\oplus\mS_i$ induces the surjective map (on edges and vertices)
 $\oplus\Ga_{\mS_i}\to\Ga_\mS$.
\eprop
The map contracts some edges of the graphs $\Ga_{\mS_i}$ (so it is not weight additive) and does not give
a significant restriction. However one has an immediate
\bcor\label{Thm Bound On r for any type}
Let $\mS\to\oplus\mS_i$ be any $\de=const$ deformation, with $r,\{r_i\}$ the corresponding numbers of branches.
Then $\bin{r}{2}\le\sum\bin{r_i}{2}$.
\ecor
This follows just from counting the number of edges in $\oplus\Ga_{\mS_i}$ and $\Ga_\mS$.

Note that this strengthens the criterion from \cite[theorem 6.1.7]{Buchweitz-Greuel80}: $r-1\le\sum(r_i-1)$.
\\\\
To get some restrictions on the possible results of deformations, start from smoothing the branches:
\bprop\label{Thm Canonical Smoothing of branches}
For every singularity type $\mS$ there exists a canonically defined type $\mS^{def}$ with smooth
branches such that $\mS$ can be deformed to $\mS^{def}$ and $\de(\mS^{def})=\de(\mS)$.
This deformation preserves the multiplicity.
\eprop
\bpr Order the branches by the minimal number of blowups needed to resolve (but the strict transform can be
tangent to the exceptional divisor). Start from the branches with the maximal such number.

Apply the minimal number of blowups till the branches become smooth
(but tangent to the exceptional divisor).
Once the branch is smooth, deform it
to intersect the divisor transversally (while the intersection numbers with other branches are preserved).
Note that this leaves all other branches intact.
Do the same with all the branches. So, we get a point with smooth branches and by construction the deformation
is canonical.
\epr
\\
So, can apply the following procedure.
\li Smoothen the branches of each singular point (canonically) in the $\de=const$ way. In this way from each graph $\Ga(C)$ we get
 $\Ga(C^{def})$ with a prescribed contraction map $\Ga(C^{def})\ra\Ga(C)$ defined as follows.

Each vertex of $\Ga(C^{def})$ corresponds to a smooth branch $C^{def}_\al$ of $C^{def}$. Under the specialization
$C^{def}\ra C$ this branch is transformed to a branch $C_\al$ of $C$. So a vertex corresponding to $C^{def}_\al$
is sent to the vertex of $C_\al$.
\\\li We have a collection of singular points with smooth branches, whose graphs have marked subgraphs. Perform
all the possible smooth-to-smooth collisions. Preserve the markings of the subgraphs.
\\\li To each resulting singularity type apply the degeneration corresponding to the contraction
of the marked subgraphs (if possible). If such a degeneration is possible and the Milnor number of
the resulting type is bigger than the sum of Milnor numbers of the initial types, such a type is
a potential candidate.
\section{Applications}\label{Sec Applications}
\subsection{Deformations of \omp s}\label{Sec Applications Deformations of OMP}
\bprop\label{Thm Application to OMP-to-OMP}
$\bullet$ Let $K_p\to(\mS_1..\mS_k)$ be a $\de=const$ deformation. Then each $\mS_i$ is an \omp.
\li
The ($\de=const$) deformation $K_p\to \cupl^k_{i=1} K_{p_i}+\Big(\bin{p}{2}-\sum\bin{p_i}{2}\Big)K_2$
with $\{p_i\ge\max(k-1,3)\}$ is possible iff $p+\bin{k}{2}\ge\sum p_i$
\eprop
\bpr $\bullet$ By proposition \ref{Thm Smooth In branches-->smooth fin branches} all the $\mS_i$ have smooth branches.
Thus the dual graph forces each $\mS_i$ to be an \omp.
\li
 Note that $\de$ of the both sides is equal.
Thus to prove sufficiency it's enough to construct in $\mP^2$ the arrangement of $p$ lines with
the prescribed combinatorics. Let $K_k$ be a (non-embedded) complete graph and
$\pi:K_k\to\mP^2$ is its projection, such that $\pi$ is injective on vertices. Let $\pi(K_k)$ be the arrangement
of lines in the plane, generated by the image of $K_k$. For a vertex $v_i\in K_k$ let $\pi(v_i)\in\pi(K_k)$
be the corresponding \omp. For each such point add $p_i-(k-1)$ lines through this point (but generic otherwise).
So one has an arrangement of $\sum\big(p_i-(k-1)\big)+\bin{k}{2}=\sum p_i-\bin{k}{2}$ lines with \omp s
of multiplicities $p_1..p_k$ and some nodes. Finally, add $p-\sum p_i+\bin{k}{2}$ generic lines. So one
has an array of $p$ lines with $\oplus K_{p_i}$ and the needed number of lines.
\\\\
\parbox{14cm}
{To prove the necessity of the condition consider the decomposition of the dual graph. We deal with \omp s, so
all the dual graphs are just the complete graphs $K_i$ (the weights of edges are 1).
Order $p_1\ge p_2\ge...$. Take a $K_{p_1}$ inside $K_p$,
now should construct the best packing of the remaining $K_{p_2}K_{p_3}...$. Any two of the embedded subgraphs
cannot intersect in more than 1 vertex. Suppose the graphs $K_{p_1}...K_{p_{i-1}}$ are embedded. So the embedded
$K_{p_i}$ can have at most $i$ vertices common with them, so $p_i-i+1$ vertices are to be added
(by the assumption $p_i\ge k-1$). Altogether this gives at least $p_1+(p_2-1)+..=\sum p_i-\bin{k}{2}$ vertices
needed to embed $K_{p_1}...K_{p_k}$, hence $p\ge\sum p_i-\bin{k}{2}$.\\\epr
}
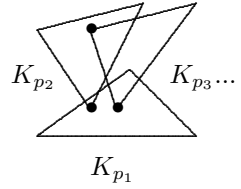
\begin{picture}(0,0)(-15,0)
\curve(0,0,60,0)\curve(0,0,35,25)\curve(35,25,60,0)

\curve(20,10,0,40)\curve(0,40,40,50)\curve(40,50,20,10)

\curve(30,10,20,40)\curve(30,10,60,50)\curve(20,40,60,50)
\put(18,8){$\bullet$}\put(28,8){$\bullet$}\put(18,38){$\bullet$}
\put(20,-15){$K_{p_1}$}\put(-10,20){$K_{p_2}$}\put(50,20){$K_{p_3}...$}

\end{picture}
\beR
Of course one can use also the semicontinuity of spectrum. If one compares the spectra of $K_p$ and $\cup K_{p_i}$
one the interval $(-\frac{1}{2}+\al,\frac{1}{2}+\al)$ and use \ref{Eq Number Of Spectra Inside Interval For OMP}
for the number of spectral pairs, one has a necessary condition:
\beq
\forall~ 0\le\al<\frac{3}{2}:~~~~
(p-1)^2-\bin{\lfloor(\frac{1}{2}-\al)p\rfloor}{2}-\bin{\lfloor(\frac{1}{2}+\al)p\rfloor}{2}\ge
\sum\Bigg((p_i-1)^2-\bin{\lfloor(\frac{1}{2}-\al)p_i\rfloor}{2}-\bin{\lfloor(\frac{1}{2}+\al)p_i\rfloor}{2}\Bigg)
\eeq
So, for example, for the deformation $K_p\to aK_3+b K_2$ this gives: $a\le \frac{2}{9}p^2-p+1$ (which is
weaker than the proposition above).
\eeR

\subsection{The canonical decomposition into \omp s}\label{Sec Applications Canonical Decomp to OMP}
\bprop
Given a type $\mS$ with smooth branches, let $\Ga_\mS=\oplus n_iK_{p_i}$ be the canonical decomposition of its
dual graph (cf. definition \ref{Def Canonical decomposition of dual graph}).

There exists a $\de=const$ deformation to the collection of \omp s: $\mS\to\bigcup n_i K_{p_i}$
(called the canonical decomposition).
The minimal number of \omp s, to which the type $\mS$
can be $\de=const$ deformed is $\sum n_i$.
\eprop
\bpr
As the canonical decomposition of the dual graph is done in steps it is enough to prove that for each step
(i.e. subtraction $\Ga_\mS-K_p$) the corresponding deformation exists.
\\
\parbox{13cm}
{Let $(C,0)$ be a representative of the type $\mS$ with the (reduced) tangent cone $T_C=(l_1..l_k)$.
Decompose the germ accordingly: $C=\bigcup C_i$, such that $T_{C_i}=l_i$. Note that each $C_i$ can be further
locally reducible.

It is enough to prove that each germ $(C_i,0)$ can be deformed into two singular points: an \omp~ at the origin
and the prescribed singularity at some other generic point (cf. the picture).
}
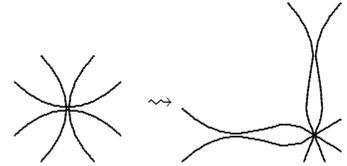
\begin{picture}(0,0)(-30,0)
\curve(-20,10,0,0.5,20,10)\curve(-20,-10,0,-0.5,20,-10)
\curve(-10,20,-0.5,0,-10,-20)\curve(10,20,0.5,0,10,-20)\put(30,0){$\rightsquigarrow$}
\end{picture}
\begin{picture}(0,0)(-90,10)
\curve(-20,10,-10,3,0,0.5,10,2,14,3,17,4,25,3,30,0,40,-6)\curve(-20,-10,-10,-3,0,-0.5,10,-2,14,-3,17,-4,25,-3,30,0,40,6)
\curve(40,50,33,40,30.5,30,32,20,32.5,16,33,10,32,5,30,0,26,-10)
\curve(20,50,27,40,29.5,30,28,20,27.5,16,27,10,28,5,30,0,34,-10)
\end{picture}
\\
So, consider one germ $C_i=\{f(x,y)=0\}$. Orient the tangent line along the $y$ axis, so that
$f=\prod_j(y(1+f_j(x,y))+x^2g_j(x))$. Here $f_j\in m_{xy}\subset\mC[[x,y]]$. Thus  $(1+f_j(x,y))$ is invertible
and the defining series of the germ can be written as $\prod_j(y+\frac{x^2g_j(x)}{1+f_j(x,y)})$. Expand in powers
of $f$, then the germ can be represented as $\prod_j(y(1+x^2f_j(x,y))+x^2g_j(x))$ (for some new $f,g$,
such that $f_j\in m_{xy}$). Iterating this procedure one arrives at the expression $\prod_j(y(1+x^{2N_j}f_j(x,y))+x^2g_j(x))$
for arbitrary large numbers $\{N_j\}$. And then, by finite determinacy, the term $x^{2N_j}f_j(x,y)$ is irrelevant.
So, can assume the germ is given in the form: $\prod_j(y+\sum_k a_{jk}x^k)$.

Consider the deformation: $f_\ep=\prod_j(y+\sum_k a_{jk}(x-\ep)^{k-1}x)$. Then at the origin $f_{\ep\ne0}$  defines
an \omp. And the germ $(f_\ep,0)$ is precisely of the type whose dual graph is $\Ga_f-K_p$.
\epr
\beR
A natural question is: whether any other  $\de=const$ deformation of a singularity into \omp s factorizes through the canonical one?
Or, at least, whether any other deformation corresponds to the further decomposition
of the dual graphs: $\bigcup n_iK_{p_i}\to...$? The following is a counterexample.

The canonical decomposition for the type $\mS=(x^4+y^{4p})$ is: $\mS\to \oplusl^p X_9$.
Suppose the deformation $\mS\to n D_4+(6p-3n)A_1$ exists. For $n\le p$ it can be factorized as
$\mS\to pX_9\to n D_4+(6p-3n)A_1$. But the case $n>p$ is the negative answer for both questions above,
since the dual graph $K_4$ (of $X_9$) does not decompose into $2K_3$ (for $D_4$).
\\
\parbox{14cm}
{It remains to show that the deformation $\mS\to n D_4+(6p-3n)A_1$ exists, e.g. for $n=p+1$.
The following construction for $n=3$, $p=2$ was given by E.Shustin.
Let a germ of curve be a line $l$ and 3 conics $C_1..C_3$ such that the conics intersect the line at three points
(so three triple points appear). The conics intersect also outside the line, adding $3A_1$.
}
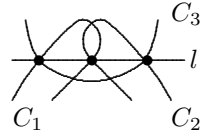
\begin{picture}(0,0)(-50,-10)
\curve(-30,0,35,0)\put(37,-3){$l$}
\curve(-30,-15,-20,0,-5,15,2,12,3,5,0,0,-5,-5,-15,-15)\put(-30,-25){$C_1$}
\curve(30,-15,20,0,5,15,-2,12,-3,5,0,0,5,-5,15,-15)\put(30,-25){$C_2$}
\curve(-25,15,-20,0,-10,-6,0,-8,10,-6,20,0,25,15)\put(30,15){$C_3$}

\put(-22.5,-2.5){$\bullet$}\put(-2.5,-2.5){$\bullet$}\put(18.5,-2.5){$\bullet$}
\end{picture}
\\\\
Note that of the four curves any pair intersects locally at two points. Thus as the three triple points
merge (and the three nodal points also join them) the family degenerates to 4 (simply) tangent curves,
i.e. the type of $x^4+y^8$.
\eeR

\beR
The proposition does not generalize to the case of the initial type $\mS$ with singular branches.
Indeed, usually there are many non-equivalent ways to smoothen the branches in a $\de=const$ way,
resulting in different dual graphs.

An interesting question is: whether each deformation of a singularity to the \omp s factorizes through
the smoothing of branches.
\eeR
\subsection{The $\de=const$ deformations of the type $x^p+y^{pk}$ to $A_k$'s}\label{Sec Applications Deformations into ADE}
Denote the type $x^p+y^{pk}$ (i.e. $p$ smooth branches, every two of them being $k-$tangent) by $K_{p,k}$.
The corresponding dual graph is the complete graph on $p$ vertices, with weights of all the edges: $k$.
\bprop\label{Thm K_pk to A_n's}
Let $K_{p,k}$ deform ($\de=const$) into a bunch of $A_i$'s. Then only $A_{2i-1}$'s appear.
Let $K_{p,k>1}\stackrel{\de=const}{\to}\bigcup n_i A_{2i-1}$. Then:
\li $n_{i>k}=0$
\li For each $i$ there exists a partition of the set $\{n_i A_{2i-1}\}$ into $\bin{p}{2}$
subsets $\{n^{(j)}_i A_{2i-1}\}_{1\le j\le\bin{p}{2}}$
such that $\suml_{j=1}^{\bin{p}{2}}n^{(j)}_i=n_i$ and $\forall j:$~ $\sum^k_{i=1}in^{(j)}_i=k$
(in particular $\sum^k_{i=1}in_i=\bin{p}{2}k$)
\li
$$\sum n_i\ge\frac{(p-1)^2k}{4}+(p-1)-\frac{1+(-1)^p}{2}~\frac{k}{4}$$.

In particular: $n_k\le\frac{(p^2-1)k}{4(k-1)}-\frac{p-1}{k-1}+\frac{1+(-1)^p}{2}~\frac{k}{4}$  and
$2n_1+n_2\ge\frac{(p-1)(p-3)k}{4}+3(p-1)-3\frac{1+(-1)^p}{2}~\frac{k}{4}$
\eprop
\bpr Consider the corresponding dual graph decomposition. Comparison of weights of edges gives $n_{i>k}=0$.
the equality $\sum^k_{i=1}in_i=\bin{p}{2}k$ is just the $\de$ of both sides.
The third inequality is obtained by the spectrum semicontinuity in the interval $(-\frac{1}{2},\frac{1}{2})$.
Comparison of the number of spectral pairs of $K_{p,k}$ (from
equation \ref{Eq Number Of Spectra Inside Interval For K_pk}) and $\bigcup n_i A_{2i-1}$ gives the result.

The bounds for $n_k,n_1,n_2$ are immediate consequence of these 3 conditions.
\epr
\beR
For specific types of the deformations above, the bounds can be slightly refined. As an example,
consider the case $p=3$. The proposition above implies: $\sum n_i\ge k+2$ and $n_k\le2$.
In fact, we have a more precise result. Consider the deformation of $K_{3,k}$ into $A_{2i-1}$'s.
As in the proposition, group them as: $\cupl_{j=1}^3\cupl_{i=1}^kn^{(j)}_i A_{2i-1}$.

\bprop
$\bullet$ $n_k\le2$ and if $n_k=2$ then the only possibility is $K_{p,k}\to 2A_{2k-1}+kA_1$
(and such deformation exists).
\li If $n_k=1$ can assume $n^{(1)}_k=1$, $n^{(2)}_k=0=n^{(3)}_k$. Let $l,m$ be the maximal integers, such that
$n^{(2)}_l\neq0\neq n^{(3)}_m$. Then $l+m\le k+1$. If $l+m=k+1$ then the deformation
$K_{p,k}\to\cup A_{2i-1}$ factorizes through $K_{p,k}\to A_{2k-1}+A_{2l-1}+A_{2m-1}+(k-1)A_1$
(and all such deformations exist).
\eprop
\bpr The additional necessary conditions are imposed by spectrum. To show the existence consider the curve
with three components:
\beq
y(y+x^k)(y+(x-t)^m(a_{l-1}x^{l-1}+..+a_0))
\eeq
For any $\{a_i\}$ this curve has $A_{2k-1}\cup A_{2m-1}$ and nodes. So, need to ensure the additional $A_{2l-1}$,
i.e. that the last two components have tangency of order $l$. Such a tangency at a given
point imposes $l+1$ conditions on $l+2$ variables ($\{a_i\}$ and $(x,y)$). Thus the system has a solution.

Finally should check the limit $t\to 0$. The limit curve has 3 components with pairwise degrees of
tangencies: $k,N,N$. Here $N\ge\max(l,m)$. If $N=k$ we have a $K_{p,k}$ point. Otherwise we have a
singularity of the topological type of $y(y+x^k)(y-x^N)$ with a bunch of nodes around. Thus the branches
can be degenerated (freely) to force the nodes to the origin and to get $K_{p,k}$.
\epr
\eeR
For low $k$ cases the equinormalizable deformations are classified below. We give only the prime
(i.e. non-factorizable deformations), all the remaining cases are obtained by further deformation.
\\
\begin{tabular}{|>{$}c<{$}|>{$}c<{$}|>{$}c<{$}|>{$}c<{$}|>{$}c<{$}|>{$}c<{$}|>{$}c<{$}|>{$}c<{$}|}
\hline
K_{3,2}:&2A_3+2A_1&A_3+4A_1&\\\hline
K_{3,3}:&2A_5+3A_1&A_5+2A_3+2A_1& 3A_3+3A_1\\\hline
K_{3,4}:&2A_7+4A_1&\ber A_7+A_5+A_3+3A_1\\ A_7+3A_3+2A_1\eer&\ber 3A_5+3A_1 \\ 2A_5+2A_3+2A_1 \\ A_5+4A_3+A_1\eer\\\hline
\end{tabular}

The classification is done by first applying the above restrictions. This leaves only the cases of the tables,
except for the candidate $K_{3,4}\to 6A_3$. This last case is ruled out by the consideration of the
deformation of the corresponding real curve.

The explicit deformations of the table are constructed starting from a real representative and then
deforming the branches.
\subsection{On the semi-continuous invariants}
The dual graph is useful in finding new semi-continuous invariants for $\de=const$ deformations.
Let $\Ga\to\oplus_i\Ga_i$ be the decomposition.
Suppose for any dual graph a function is defined $f:\Ga\to\mZ$. Then can compare $f(\Ga)$ vs $\sum f(\Ga_i)$.

\bex
$\bullet$ By counting the number of branches one has: $\bin{r}{2}\le\sum\bin{r_i}{2}$.
\li Let $w_i$ be the weights of the graph, choose $f:=(\sum w_i^p)^{\frac{1}{p}}$. Then Minkowski's inequality gives:
$f(\Ga)\le\sum f(\Ga_i)$.

\li Let $f:=\sum w_i^p$. Then $f(\Ga)\le r^{p-1}\sum f(\Ga_i)$.
\eex

\end{document}